\documentclass{article} 
\usepackage[nonatbib,preprint]{helpers/neurips_2023}
\usepackage[backend=biber,style=numeric,sorting=nty]{biblatex}
\bibliography{biblio.bib}

\usepackage[utf8]{inputenc} 
\usepackage[T1]{fontenc}    
\usepackage{hyperref}       
\usepackage{url}            
\usepackage{booktabs}       
\usepackage{amsfonts}       
\usepackage{nicefrac}       
\usepackage{microtype}      
\usepackage{xcolor}         

\usepackage[algo2e,noend,boxruled]{algorithm2e}
\usepackage{algorithm}
\usepackage{amsmath}
\usepackage{amssymb}
\usepackage{mathtools}
\usepackage{amsthm}
\usepackage{amsfonts}
\usepackage{amsopn}
\usepackage{array}
\usepackage{bm}
\usepackage[cal=boondoxo]{mathalfa}
\usepackage{multirow}
\usepackage{ragged2e}
\usepackage{tabularx}
\usepackage{ulem}
\usepackage{verbatim}
\usepackage{mathtools}
\mathtoolsset{showonlyrefs}

\newcommand{\fat}[1]{\ifmmode\bm{#1}\else\textbf{#1}\fi}

\newcommand{\set}[1]{\mathbb{#1}}

\newcommand{\vect}[1]{\ensuremath{\fat{#1}}}

\newcommand{\matr}[1]{#1}

\newcommand{\tens}[1]{\mathcal{#1}}

\newcommand{\func}[1]{\textsf{#1}}

\newcommand{\ffunc}[1]{\widehat{#1}}

\newcommand{\order}[1]{\mathcal{O}\left( #1 \right)}
\newcommand{\orders}[2]{\mathcal{O}\csname#1l\endcsname( #2 \csname#1r\endcsname)}

\newcommand{\ff}[0]{\func{f}}

\newcommand{\fl}[0]{\ffunc{L}}

\newcommand{\ty}[0]{\tens{Y}}

\newcommand{\distp}[0]{\func{p}}

\newcommand{\tp}[0]{\tens{P}}

\newcommand{\tg}[0]{\tens{G}}

\newcommand{\vx}[0]{\vect{x}}

\newcommand{\vz}[0]{\vect{z}}

\newcommand{\vxmin}[0]{\vx_{min}}

\newcommand{\ymin}[0]{y_{min}}

\newcommand{\neff}[0]{\overline{N}}

\newcommand{\reff}[0]{\overline{R}}

\newcommand{\trans}[1]{\textsf{F}[#1]}

\newcommand{\ocf}[0]{\func{g}}

\newcommand{\oct}[0]{\func{G}}

\def\ie{i.\,e.}
\title{
    PROTES: Probabilistic Optimization with Tensor Sampling
}
\author{
    Anastasia Batsheva$^*$
        \\
        Skolkovo Institute of Science and Technology
        \\
        Moscow, Russia
        \\
        \texttt{a.batsheva@skoltech.ru}
    \And
    Andrei Chertkov$^*$
        \\
        Skolkovo Institute of Science and Technology
        \\
        Moscow, Russia
        \\
        \texttt{a.chertkov@skoltech.ru}
    \AND
    Gleb Ryzhakov$^*$
        \\
        Skolkovo Institute of Science and Technology
        \\
        Moscow, Russia
        \\
        \texttt{g.ryzhakov@skoltech.ru}
    \And
    Ivan Oseledets
        \\
        Skolkovo Institute of Science and Technology
        \\
        and AIRI, Moscow, Russia
        \\
        \texttt{i.oseledets@skoltech.ru}
}
\begin{document}
\maketitle
\begin{abstract}
    We developed a new method PROTES for black-box optimization, which is based on the probabilistic sampling from a probability density function given in the low-parametric tensor train format.
    We tested it on complex multidimensional arrays and discretized multivariable functions taken, among others, from real-world applications, including unconstrained binary optimization and optimal control problems, for which the possible number of elements is up to $2^{100}$.
    In numerical experiments, both on analytic model functions and on complex problems, PROTES outperforms existing popular discrete optimization methods (Particle Swarm Optimization, Covariance Matrix Adaptation, Differential Evolution, and others).
\end{abstract}
\def\thefootnote{*}\footnotetext{Equal contribution.}\def\thefootnote{\arabic{footnote}}
\section{Introduction}
    \label{sec:introduction}
    The multidimensional optimization problem is one of the most common in machine learning.
It includes the important case of gradient-free discrete optimization~\cite{parker2014discrete, lauand2022approaching, lin2022gradientfree, wang2022zarts}:
\begin{equation}\label{eq:task}
    \vxmin = 
        \min\limits_{\vx} \ff(\vx),
    \quad
    \textit{ s.t. }
    \vx = [n_1, n_2, \ldots, n_d],
    \quad
    n_i \in \{ 1, 2, \ldots, N_i \},
\end{equation}
where $d$ is the dimensionality of the problem, and $N_1, N_2, \ldots, N_d$ are the numbers of items for each dimension.
Such settings arise when searching for the minimum or maximum element in an implicitly given multidimensional array (tensor\footnote{
    A tensor is just a multidimensional array with several dimensions $d$ ($d \geq 1$).
    A two-dimensional tensor ($d = 2$) is a matrix, and when $d = 1$ it is a vector.
    For scalars we use normal font ($a, b, c, \ldots$), we denote vectors with bold letters ($\vect{a}, \vect{b}, \vect{c}, \ldots$), we use upper case letters ($\matr{A}, \matr{B}, \matr{C}, \ldots$) for matrices, and calligraphic upper case letters ($\tens{A}, \tens{B}, \tens{C}, \ldots$) for tensors with $d > 2$.
    The $(n_1, n_2, \ldots, n_d)$th entry of a $d$-dimensional tensor $\ty \in \set{R}^{N_1 \times N_2 \times \ldots \times N_d}$ is denoted by $y = \ty[n_1, n_2, \ldots, n_d]$, where $n_i = 1, 2, \ldots, N_i$ ($i = 1, 2, \ldots, d$), and~$N_i$ is a size of the $i$-th mode.
    The mode-$i$ slice of such tensor is denoted by $\ty[n_1, \ldots, n_{i-1}, :, n_{i+1}, \ldots, n_d]$.
}), including when considering the discretization of functions from a continuous argument.
Multidimensional discrete optimization problems are still computationally difficult in the case of complex target functions or large dimensions, and efficient direct gradient-free optimization procedures are highly needed.

The development of methods for low-rank tensor approximations has made it possible to introduce fundamentally new approaches for the approximation, storage, and operation with multidimensional tensors~\cite{grasedyck2013literature, cichocki2016tensor, cichocki2017tensor, wang2023tensor}.
One of the common methods for low-rank approximation is the tensor train (TT) decomposition~\cite{oseledets2011tensor}, which allows bypassing the curse of dimensionality.
Many useful algorithms (e.\,g., element-wise addition, multiplication, solution of linear systems, convolution, 
integration, etc.) have effective implementations for tensors given in the TT-format.
The complexity of these algorithms turns out to be polynomial in the dimension $d$ and the mode size $N$.
It makes the TT-decomposition extremely popular in a wide range of applications, including computational mathematics~\cite{parcollet2023learning, alexandrov2023challenging} and machine learning~\cite{qi2023exploiting, kour2023efficient}.

In the last few years, new discrete optimization algorithms based on the TT-format have been proposed: TTOpt~\cite{sozykin2022ttopt}, Optima-TT~\cite{chertkov2022optimization}, and several modifications~\cite{selvanayagam2022global, shetty2022tensor, nikitin2022quantum, soley2023global}.
However, the development of new more accurate, and robust TT-based methods for multidimensional discrete optimization is possible.
In this work, we extend recent approaches for working with probability distributions and sampling in the TT-format~\cite{dolgov2020approximation, novikov2021tensor} to the optimization area.
That is, we specify a multidimensional discrete probability distribution in the TT-format, followed by efficient sampling from it and updating its parameters to approximate the optimum in a better way.
This makes it possible to build an effective optimization method, and the contributions of our work are as follows:
\begin{itemize}
    \item
        We develop a new method PROTES for optimization (finding the minimum or maximum\footnote{
            Further, for concreteness, we will consider the minimization problem in this paper, while the proposed method can be applied to the discrete maximization problem without any significant modifications.
        } value) of multidimensional data arrays and discretized multivariable functions based on a sampling method from a probability distribution in the TT-format;
    \item
        We apply\footnote{
            The program code with the proposed method and numerical examples is available in the public repository \url{https://github.com/anabatsh/PROTES}.
        } PROTES for various analytic model functions and for several multidimensional QUBO and optimal control problems to demonstrate its performance and compare it with popular discrete optimization algorithms (Particle Swarm Optimization, Covariance Matrix Adaptation, Differential Evolution and NoisyBandit) as well as TT-based methods (TTOpt and Optima-TT).
        We used the same set of hyperparameters of our algorithm for all experiments, and obtained the best result for $19$ of the $20$ problems considered.
\end{itemize}
\section{Optimization with probabilistic sampling}
    \label{sec:method}
    Our problem is to minimize the given multivariable discrete black-box function $\ff$~\eqref{eq:task}.
It can be formulated in terms of the multi-index search in an implicitly given $d$-dimensional tensor
\begin{equation}\label{eq:task_discrete}
    \ty \in \set{R}^{
        N_1 \times N_2 \times \ldots \times N_d},
    \quad
    \ty[n_1,\, n_2,\, \ldots,\, n_d] =
        \ff(\vx),
    \quad
    \vx = [n_1,\, n_2,\, \ldots,\, n_d],
\end{equation}
for all $n_i = 1, 2, \ldots, N_i$ ($i = 1, 2, \ldots, d$).
The essence of our idea is to use a probabilistic method to find the minimum $\vxmin$.
We propose to establish a discrete distribution $\distp(\vx)$ from which the minimum could be sampled with high probability.
This distribution can be specified as a tensor~$\tp_\theta \in \set{R}^{N_1 \times N_2 \times \ldots \times N_d}$ in some low-parametric representation with a set of parameters~$\theta$, having the same shape as the target tensor $\ty$. 

\begin{figure}[t!]
    \centering
    \includegraphics
        [width=0.99\linewidth]
        {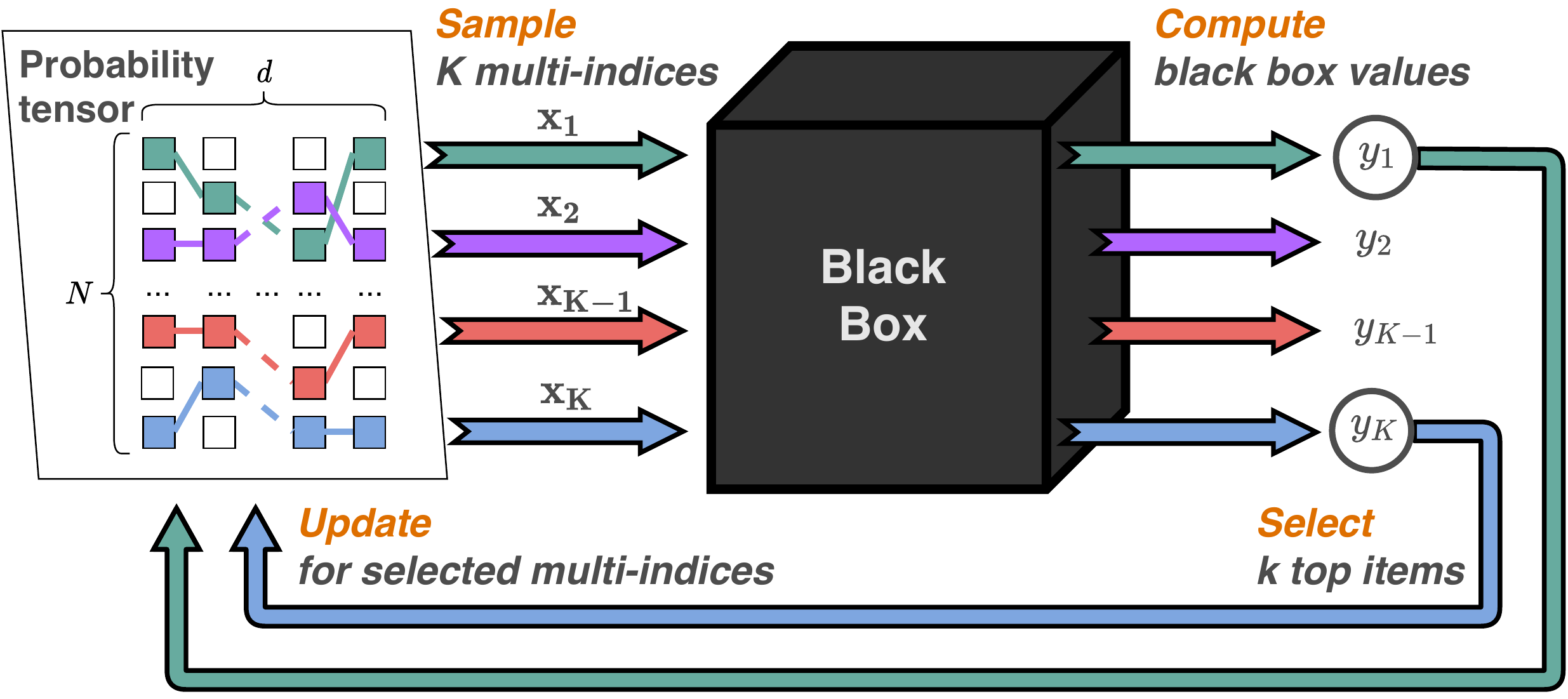}
    \caption{
        Schematic representation of the proposed optimization method PROTES.
    }
    \label{fig:protes}
\end{figure}

We start from a random non-negative tensor $\tp_\theta$ 
and iteratively perform the following steps 
until the budget is exhausted or until convergence 
(see graphic illustration in Figure~\ref{fig:protes}):
\begin{enumerate}
    \item \textbf{Sample}
        $K$ candidates of $\vxmin$ from the current distribution $\tp_\theta$:
        $\mathcal{X}_K = \{\vx_1, \vx_2, \ldots, \vx_K\}$;
    \item \textbf{Compute}
        the corresponding function values:
        $y_1 = \ff(\vx_1)$, $y_2 = \ff(\vx_2)$, \ldots, $y_K = \ff(\vx_K)$;
    \item \textbf{Select}
        $k$ best candidates with indices~$\mathcal{S} = \{ s_1, s_2, \ldots, s_k \}$ 
        from $\mathcal{X}_K$ with the minimal objective value, \ie, $y_{j} \leq y_{J}$ for all $j \in \mathcal{S}$ and $J \in \{ 1,\, 2,\, \ldots,\, K \} \setminus \mathcal{S}$;
    \item \textbf{Update}
        the probability distribution $\tp_\theta$ ($\theta \rightarrow \theta^{(new)}$) to increase the likelihood of selected candidates $\mathcal{S}$.
        We make several ($k_{gd}$) gradient ascent steps with the learning rate $\lambda$ for the tensor $\tp_\theta$, using the following loss function
        \begin{equation}\label{eq:loss-function}
            \fl_\theta(\{x_{s_1},\,x_{s_2},\,\ldots,\,x_{s_k}\}) = \sum_{i=1}^k
                \log \left(
                    \tp_\theta[\vx_{s_i}]
                \right).
        \end{equation}
\end{enumerate}
After a sufficient number of iterations, 
we expect the tensor $\tp_\theta$ to represent an almost 
Kronecker delta-function with a pronounced peak in the value of the minimum of the target function 
(or several peaks if the minimum is not unique).
Therefore, this value will be sampled during the steps of our algorithm since the probability of sampling other values will be sufficiently small.

From a low-parameter representation $\tp_\theta$ we expect an efficient sampling algorithm and efficient calculation procedure for the logarithms in~\eqref{eq:loss-function} with automatic differentiation capability to enable gradient ascent methods. 
As will be shown below, the TT-representation of tensors satisfies these requirements.
Further, we will omit the index~$\theta$, assuming that the parameterized representation of the tensor~$\tp$ corresponds to the TT-format.
Note that the values~$K$, $k$, $k_{gd}$, $\lambda$ and the number of parameters in~$\theta$ (\ie, rank of the TT-decomposition) are the hyperparameters of our algorithm.

\section{Basic properties of the tensor train format}
    \label{sec:tensortrain}
\begin{figure}[t!]
    \centering
    \includegraphics
        [width=0.99\linewidth]
        {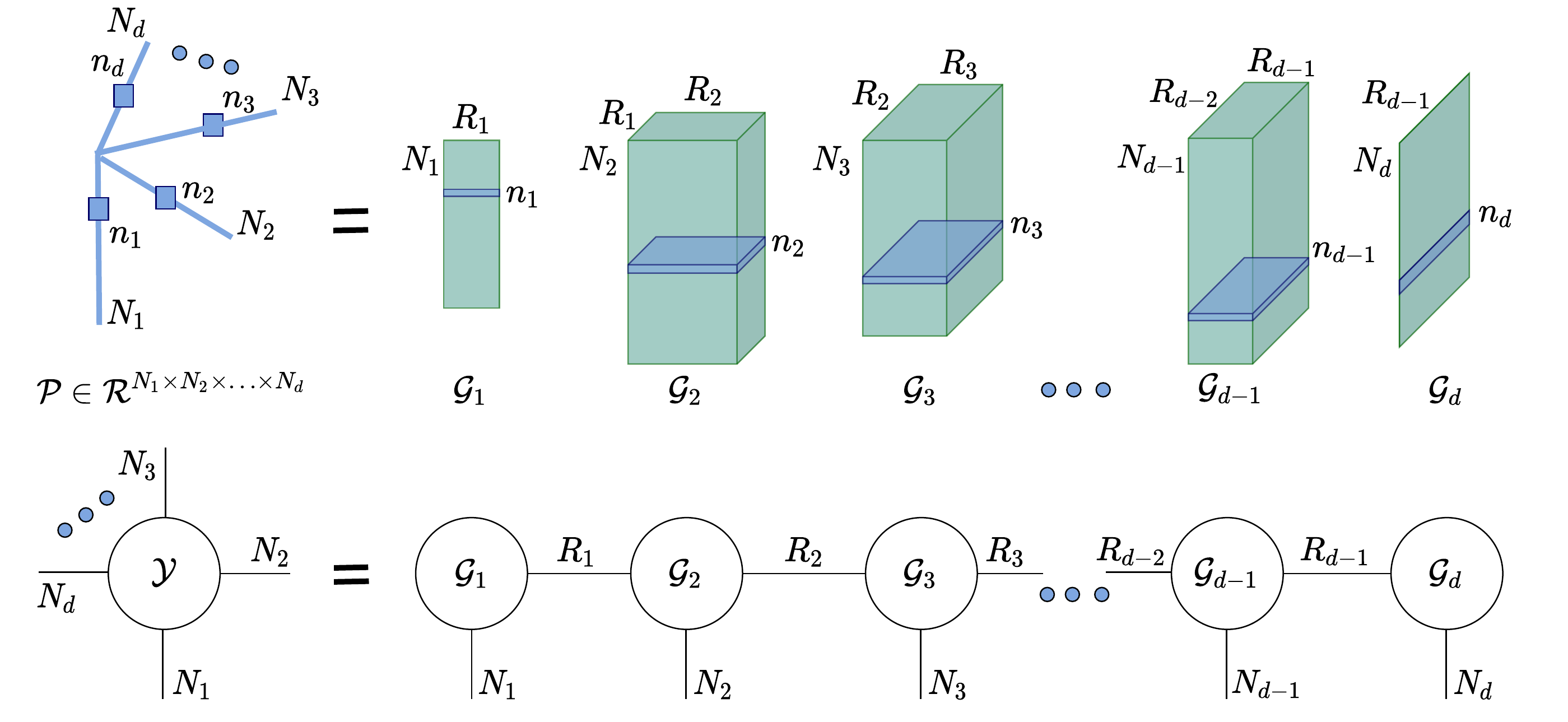}
    \caption{
        Schematic representation of the TT-decomposition.
        The top picture demonstrates the calculation of the specific tensor element $\vx = [n_1, n_2, \ldots, n_d]$ from its TT-representation, and the bottom picture presents the related tensor network diagram.
    }
    \label{fig:tt-element}
\end{figure}

    Let us dwell on the concept of the TT-format.
A $d$-dimensional tensor $\tp \in \set{R}^{N_1 \times N_2 \times \ldots \times N_d}$ is said to be in the TT-format~\cite{oseledets2011tensor}
if its elements are represented by the following formula
\begin{equation}\label{eq:tt-repr-tns}
    \tp [n_1, n_2, \ldots, n_d]
    =
    \sum_{r_1=1}^{R_1}
    \sum_{r_2=1}^{R_2}
    \cdots
    \sum_{r_{d-1}=1}^{R_{d-1}}
        \tg_1 [1, n_1, r_1]
        \;
        \tg_2 [r_1, n_2, r_2]
        \;
        \ldots
        \;
        \tg_d [r_{d-1}, n_d, 1],
\end{equation}
where $(n_1, n_2, \ldots, n_d)$ is a multi-index ($n_i = 1, 2, \ldots, N_i$ for $i = 1, 2, \ldots, d$), integers $R_{0}, R_{1}, \ldots, R_{d}$ (with convention $R_{0} = R_{d} = 1$) are named TT-ranks, and three-dimensional tensors $\tg_i \in \set{R}^{R_{i-1} \times N_i \times R_i}$ ($i = 1, 2, \ldots, d$) are named TT-cores.
The TT-decomposition~\eqref{eq:tt-repr-tns} (see also an illustration in Figure~\ref{fig:tt-element}) allows to represent a tensor or a discretized multivariable function in a compact and descriptive low-parameter form, which is linear in dimension~$d$, \ie, it has less than $d \cdot \max_{i=1,\ldots,d}(N_i R_i^2) \sim d \cdot \neff \cdot \reff^2$ parameters, where $\neff$ and $\reff$ are effective (``average'') mode size and TT-rank respectively.

Linear algebra operations (e.\,g., element-wise addition, solution of linear systems, convolution, integration, etc.) on tensors in the TT-format, respectively, also have complexity linear in dimension if the TT-ranks are bounded.
The TT-approximation for a tensor or discretized multivariable function may be built by efficient numerical methods, e.\,g., TT-SVD~\cite{oseledets2009breaking}, TT-ALS~\cite{chertkov2022black}, and TT-cross~\cite{oseledets2010ttcross}.
A~detailed description of the TT-format and related algorithms are given in works~\cite{oseledets2011tensor, cichocki2016tensor, cichocki2017tensor}.
Below we discuss only three operations in the TT-format, which will be used later in the work.

\paragraph{Construction of the random TT-tensor.}
To build a random non-negative TT-tensor of a given size $(N_1, N_2, \ldots, N_d)$ with a constant TT-rank $R$, it is enough to generate $d$ TT-cores $\tg_1, \tg_2, \ldots, \tg_d$ ($3$-dimensional tensors) with random elements from the uniform distribution on the interval $(0, 1)$.
We will refer to this method as $\func{tt\_random}(R, \,[N_1, N_2, \ldots, N_d])$.

\paragraph{Computation of the log-likelihood in the TT-format.}
To calculate the logarithm $\log \tp[\vx]$ in a given multi-index $\vx = (n_1, n_2, \ldots, n_d)$, we can use the basic formula~\eqref{eq:tt-repr-tns} and then take the logarithm of the result.
It can be shown that this operation has complexity $\orders{Big}{d \cdot \reff^2}$,
because, roughly speaking, we $(d-1)$ times multiply a vector of length~$\reff$ by a matrix of size~$\reff \times \reff$  to get the result.
The corresponding method will be called below as $\func{tt\_log}(\tp, \vx)$.

\paragraph{Sampling from the tensor in the TT-format.}
To generate a multi-index $\vx$ with a probability proportional to the corresponding value $p = \tp[\vx]$ of the TT-tensor $\tp$, we use the approach proposed in the work~\cite{dolgov2020approximation}.
The method is based on the sequential calculation of univariate conditional densities with efficient integration in the TT-format, and the estimate for its complexity turns out to be the following:
$
    \orders{Big}{
        K \cdot d \cdot (\neff + \reff) \cdot \reff
        + K\cdot d\cdot\alpha(\neff)
        }
$,
where $K$ is a number of requested samples, and $\alpha(n)$ is a complexity of sampling from generalized Bernoulli distribution with~$n$ outcomes.
Note that the algorithm allows sampling in the case of the initially non-normalized tensor, so we don't have to calculate the normalization factor.
We will refer to this method as $\func{tt\_sample}(\tp, K)$.
\section{Optimization method PROTES}
    \label{sec:scheme}
    With the formal scheme of the proposed approach given in Section~\ref{sec:method} and the description of operations \func{tt\_random}, \func{tt\_log} and \func{tt\_sample} given in Section~\ref{sec:tensortrain}, we can formulate our method PROTES for gradient-free discrete optimization in the TT-format, as presented in Algorithm~\ref{alg:protes}.
We denote as \func{adam}, a function that performs $k_{gd}$ steps of gradient ascent for the TT-tensor $\tp$ at multi-indices~$\mathcal{X}$ by the well-known Adam method~\cite{kingma2014adam}.
In this case, the learning rate is $\lambda$, the loss function is given in~\eqref{eq:loss-function}, and \func{tt\_log} with automatic differentiation support is used for the log-likelihood computation.

\begin{algorithm}[t!]
\small
\SetAlgoLined
\LinesNumbered
\KwData{
    the function $\ff(\vx)$, that computes the value of the target tensor $\ty \in \set{R}^{N_1 \times N_2 \times \ldots \times N_d}$ at the multi-index $\vx = [n_1, n_2, \ldots, n_d]$;
    the maximum number of requests $M$;
    the number of generated samples per iteration $K$;
    the number of selected candidates per iteration $k$;
    the number of gradient ascent steps $k_{gd}$;
    the gradient ascent learning rate $\lambda$;
    the TT-rank of the probability tensor $R$.
}
\KwResult{
    $d$-dimensional multi-index $\vxmin$, which relates to
    the minimum value of the tensor $\ty$.
}

Initialize target multi-index and tensor value:
$
    \vxmin = \texttt{None},
    \,
    \ymin = \infty
$

Generate random non-negative rank-$R$ TT-tensor:
$
    \tp = \func{tt\_random}\left(
        R, \,[N_1, N_2, \ldots, N_d]
    \right)
$

\For{$iter = 1$ \KwTo $M/K$}{
    Generate $K$ samples from $\tp$:
    $
        \vx_1, \vx_2, \ldots, \vx_K =
            \func{tt\_sample}(\tp, K)
    $

    Compute related tensor values:
    $
        y_1 = \ff(\vx_1), 
        y_2 = \ff(\vx_2), 
        \ldots,
        y_K = \ff(\vx_K) 
    $

    Find indices
    $
        \mathcal{S} = \{
            s_1, s_2, \ldots, s_k
        \}
    $
    for top-$k$ minimum items in the list
    $
        [y_1, y_2, \ldots, y_K]
    $

    Collect k-top candidates:
    $
        \mathcal{X} = \{
            \vx_{s_1},
            \vx_{s_2},
            \ldots,
            \vx_{s_k}
        \},
        \,\,
        \mathcal{Y} = \{
            y_{s_1},
            y_{s_2},
            \ldots,
            y_{s_k}
        \}
    $
    
    If $\mathcal{Y}$ contains a value less than $\ymin$, then update $\vxmin$ and $\ymin$

    Update the TT-tensor:
    $
        \tp \gets \func{adam}(
            \tp, \fl, \mathcal{X}, k_{gd}, \lambda)
    $
    // \textit{with the loss function~\eqref{eq:loss-function} and the method \func{tt\_log}}
}

\Return{$\vxmin$}.

\caption{Method PROTES in the TT-format for multidimensional discrete black-box optimization}
\label{alg:protes}
\end{algorithm}

\paragraph{Computational complexity of the method.}
Let estimate the complexity of the proposed algorithm, assuming that the number of requests to the target function (black-box) $M$ is fixed.
With the known estimate for the complexity of the \func{tt\_sample} function, we can obtain the complexity of the sampling operations:
$
    \order{
        \frac MK \cdot K \cdot d \cdot \left(
            (\neff + \reff) \cdot \reff
            +
            \alpha(\neff)
        \right)
    }
$.
Assuming that the complexity of one gradient step coincides with the complexity of calculating the differentiated function and using the estimate for the \func{tt\_sample} function, we can estimate the total complexity of the tensor updates:
$
    \order{
        \frac MK\cdot k \cdot k_{gd} \cdot
        d \cdot\reff^2
    }
$.
Combining the two above estimates we obtain the complexity of the method
\begin{equation}\label{eq:complexity}
    \order{
        M \cdot d\cdot\Bigl(\frac kK \cdot k_{gd} \cdot
        \reff^2
        +(\neff+\reff)\cdot\reff+\alpha(N)
        \Bigr)
    }.
\end{equation}
However, it is important to note that this estimate does not take into account the complexity of calculating $M$ times the objective function $\ff$, which in practical applications can be significant and many times greater than the estimate~\eqref{eq:complexity}.

\paragraph{The intuition behind the method.}
The proposed method PROTES, like most gradient-free optimization approaches, is empirical, however, we can establish its connection with a well-known REINFORCE trick algorithm~\cite{williams1992simple}.
Let make a monotonic transformation~$\trans{\ff}(\vx)$ of the target function~$\ff$ to be minimized that transforms minimum to maximum.
A reasonable choice for $\trans{\cdot}$ is the Fermi-Dirac function
\begin{equation}\label{eq:fermi_dirac}
    \trans{\ff}(\vx) = \frac1{
        \exp \bigl(
            (\ff(\vx) - y_{\text{min}} - E) / T
        \bigr)
        + 1
    },
\end{equation}
where $y_{\text{min}}$ is an exact or approximate minimum of $\ff$, $T>0$ is a parameter and $E>0$ is some small threshold.
With the function~$\trans{\ff}$ we can find a maximum of the expectation
$
    \max_\theta
        \mathsf E_{\xi_\theta} \trans{\ff}(\xi_\theta)
$,
where a family of random variables~$\xi_\theta$ has a parameterised distribution function~$\func{p}_\theta(\vx)$.
Using REINFORCE trick, we can estimate the gradient of the expectation by the following Monte-Carlo-like expression
\begin{equation}\label{eq:RF_t:MC}
    \nabla_\theta
    \mathsf E_{\xi_\theta} \trans{\ff}(\xi_\theta) 
    \approx
    \frac1M\sum_{i=1}^M\trans{\ff}(\vx_i)\nabla_\theta \log \func{p}_\theta(\vx_i),
\end{equation}
where $\{\vx_i\}_1^M$ are independent realizations of the random variable~$\xi_\theta$.
If we find the optimal values of~$\theta$, then we expect the optimal distribution~$\func{p}_\theta$ to have a peak at the point of maximum for function~$\trans{\ff}$.
Thus we can obtain the argument of its maximum by sampling from this distribution.
For very small values of~$T$, only a few terms contribute to the sum~\eqref{eq:RF_t:MC}, namely those~$\vx_i$ for which $\ff(\vx_i) - y_{min} < E$ is hold.
For these values of~$\vx$, $\trans{\ff}$ is close to~$1$, while for the other samples its value is~$0$.
Hence, we can discard all other samples and keep a few samples with the best values.
So, we come to the loss function~\eqref{eq:loss-function}, where instead of the parameter~$E$ we use a fixed number~$k$ of the best samples, \ie, the samples for which the value of the target function~$\ff$ is the smallest.

\paragraph{Application of the method to constrained optimization.}
A very nice property of the proposed method is that it can be adapted to efficiently handle constraints such as a specified set of admissible multi-indices.
One option is just to remove invalid samples from the top-k values, but in some cases, the probability of sampling multi-indices that are admissible is very low, so this approach will not work.
Instead, if the constraint permits, we use the algorithm from the work~\cite{ryzhakov2023constructive} for the constructive building of tensors in the TT-format by a known analytic function, which defines the constraints.
Once the indicator tensor ($1$ if the index is admissible and $0$ if it is not) is built in the TT-format, we can just initialize the starting distribution $\tp$ by it, and it will be guaranteed that the samples almost always belong to the admissible set.
\section{Related work}
    \label{sec:related}
    Below we give a brief analysis of classical approaches for discrete optimization and then discuss the methods based on the low-rank tensor approximations, which have become popular in the last years.

\paragraph{Classical methods for gradient-free optimization.}
In many situations, the problem-specific target function is not differentiable, too complex, or its gradients are not helpful due to the non-convex nature of the problem~\cite{alarie2021two}, and standard well-known gradient-based methods cannot be applied directly.
The examples include hyper-parameter selection, training neural networks with discrete weights, and policy optimization in reinforcement learning.
In all these contexts, efficient direct gradient-free optimization procedures are highly needed.
In the case of high dimensional black-box optimization, evolutionary strategies (ES)~\cite{doerr2021survey} are one of the most advanced methods.
This approach aims to optimize the parameters of the search distribution, typically a multidimensional Gaussian, to maximize the objective function.
Finite difference schemes are commonly used to approximate gradients of the search distribution.
Numerous works proposed techniques to improve the convergence of ES~\cite{nesterov2017random}, for example, second-order natural gradient~\cite{wierstra2014natural} or the history of recent updates (Covariance Matrix Adaptation Evolution Strategy; CMA-ES)~\cite{hansen2006cma} may be used to generate updates.
There is also a large variety of other heuristic methods for finding the global extremum.
In particular, we note such popular approaches as NoisyBandit~\cite{scarlett2017lower}, Particle Swarm Optimization (PSO)~\cite{kennedy1995particle}, Simultaneous Perturbation Stochastic Approximation (SPSA)~\cite{maryak2001global}, Differential Evolution (DE)~\cite{storn1997differential} and scrambled-Hammersley (scr-Hammersley)~\cite{hammersley1960monte}.

\paragraph{Tensor-based methods for gradient-free optimization.}
Recently, the TT-decomposition has been actively used for multidimensional optimization.
An iterative method TTOpt based on the maximum volume approach is proposed in the work~\cite{sozykin2022ttopt}.
TTOpt utilizes the theorem of sufficient proximity of the maximum modulo element of the submatrix having the maximum modulus of the determinant to the maximum modulo element of the tensor.
Based on this observation, tensor elements are sampled from specially selected successive unfoldings of the tensor.
To be able to find the minimum element, dynamic mapping of the tensor elements is carried out, which converts the minimum values into maximum ones.
The authors applied this approach to the problem of optimizing the weights of neural networks in the framework of reinforcement learning problems in~\cite{sozykin2022ttopt} and to the QUBO problem in~\cite{nikitin2022quantum}.
A similar optimization approach was also considered in~\cite{selvanayagam2022global} and~\cite{shetty2022tensor}.
One more promising algorithm, named Optima-TT, was proposed in recent work~\cite{chertkov2022optimization}.
This approach is based on the probabilistic sampling from the TT-tensor and makes it possible to obtain a very accurate approximation for the optimum of the given TT-tensor.
However, this method is intended for directly optimizing the TT-tensors, which means that its success strongly depends on the quality of the TT-approximation for the original multidimensional data array.
Therefore, one of the related methods in the TT-format (TT-SVD, TT-ALS, TT-cross, etc.) should be additionally used for approximation.
\section{Numerical experiments}
    \label{sec:experiments}
\begin{figure}[t!]
    \centering
    \includegraphics[width=0.99\linewidth]
        {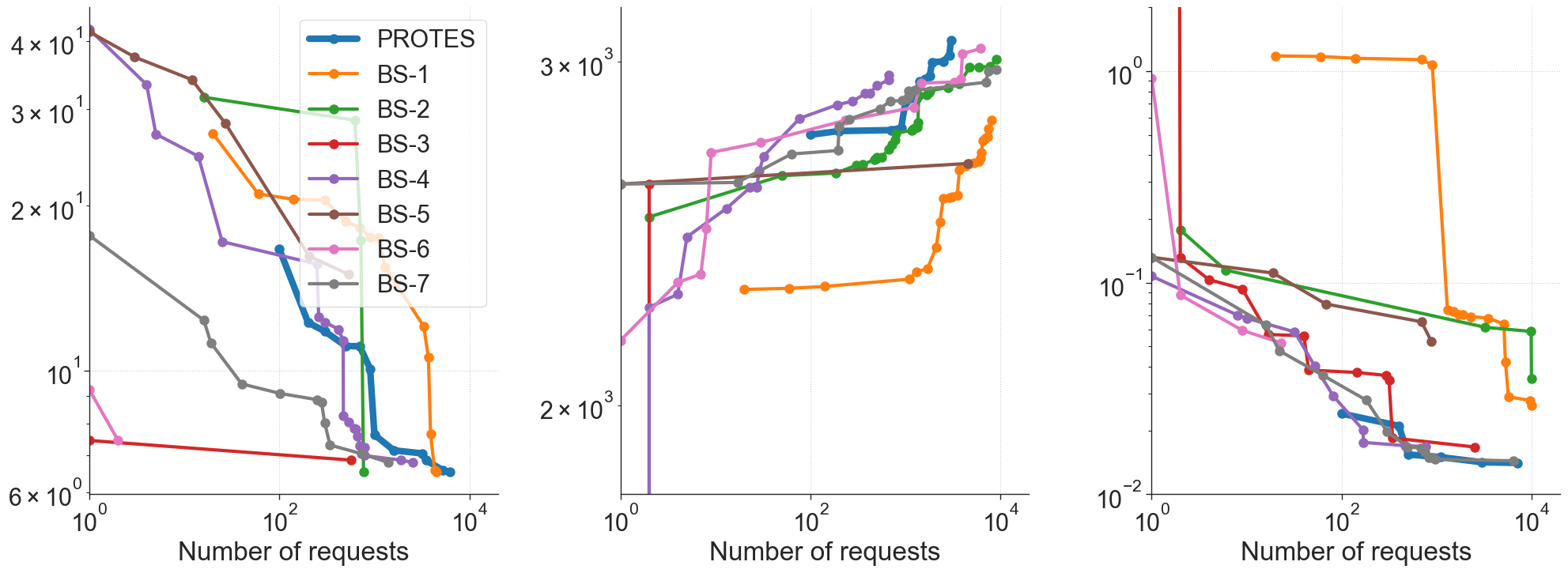}
    \caption{
        The dependence of the found optimum value $\ff(\vxmin)$ in log scale on the number of requests to the objective function for all considered optimization methods for benchmarks P-02 (plot on the left), P-14 (plot in the middle; for clarity of demonstration, the negative values for this benchmark are inverted, \ie, we plot $-\ff(\vxmin)$ values) and P-16 (plot on the right).
    }
    \label{fig:deps}
\end{figure}

    {
    To evaluate the effectiveness of the proposed method, we carried out a series of $20$ numerical experiments for various formulations of model problems.
    The results are presented in Table~\ref{tbl:results}, where we report the approximation to the minimum value for each model problem (P-1, P-2, \ldots, P-20) and all considered optimization methods (PROTES, BS-1, BS-2, \ldots, BS-7).
    Taking into account the analysis of discrete optimization methods in the previous section, as baselines we consider two tensor-based optimization methods: TTOpt\footnote{
        We used implementation of the method from \url{https://github.com/AndreiChertkov/ttopt}.
    } (BS1) and Optima-TT\footnote{
        We used implementation from \url{https://github.com/AndreiChertkov/teneva}.
        The TT-tensor for optimization was generated by the TT-cross method.
    }  (BS2), and five popular gradient-free optimization algorithms from the nevergrad framework:~\cite{bennet2021nevergrad}\footnote{
        See \url{https://github.com/facebookresearch/nevergrad}.
    } OnePlusOne (BS3), PSO (BS4), NoisyBandit (BS5), SPSA (BS6), and Portfolio approach (BS7), which is based on the combination of CMA-ES, DE, and scr-Hammersley methods.
    The model problems and obtained results will be discussed in detail below in this section.
    
    For all the considered optimization problems, we used the default set of parameters for baselines, and for PROTES we fixed parameters as $K=100$, $k=10$, $k_{gd}=1$, $\lambda=0.05$, $R=5$ (the description of these parameters was presented in Algorithm~\ref{alg:protes}).
    For all methods, the limit on the number of requests to the objective function was fixed at the value  $M=10^4$.
    As can be seen from Table~\ref{tbl:results}, PROTES, in contrast to alternative approaches, gives a consistently top result for almost all model problems (the best result for $19$ of the $20$ problems considered).
    To demonstrate the convergence behavior of methods, we also present the corresponding plot in Figure~\ref{fig:deps}.
}

\begin{table}[t!]
\caption{
    Minimization results for all selected benchmarks (P-01 -- P-20).
    The values obtained by the proposed method PROTES and by all considered baselines (BS1 -- BS7) are reported.
}
\label{tbl:results}
\renewcommand{\arraystretch}{1.3}
\begin{center}
\begin{tiny}
\begin{sc}
\begin{tabular}
{p{1.5cm}p{0.7cm}rrrrrrrr}\hline

        &
        &
PROTES  & 
BS-1    & 
BS-2    &
BS-3    &
BS-4    &
BS-5    &
BS-6    &
BS-7    \\ \hline

\multirow{10}{*}{\parbox{1.3cm}{Analytic Functions}}
    & P-01
        &  \fat{1.3e+01}
        &  \fat{1.3e+01}
        &  \fat{1.3e+01}
        &  \fat{1.3e+01}
        &  \fat{1.3e+01}
        &  2.1e+01
        &  \fat{1.3e+01}
        &  \fat{1.3e+01}
    \\ 
    & P-02
        &  \fat{6.5e+00}
        &  \fat{6.5e+00}
        &  \fat{6.5e+00}
        &  6.9e+00
        &  6.8e+00
        &  1.5e+01
        &  7.5e+00
        &  6.8e+00
    \\ 
    & P-03
        & \fat{-9.4e-01}
        & \fat{-9.4e-01}
        & \fat{-9.4e-01}
        & \fat{-9.4e-01}
        & \fat{-9.4e-01}
        & -3.5e-01
        & \fat{-9.4e-01}
        & \fat{-9.4e-01}
    \\ 
    & P-04
        &  \fat{1.3e+00}
        &  \fat{1.3e+00}
        &  \fat{1.3e+00}
        &  \fat{1.3e+00}
        &  \fat{1.3e+00}
        &  6.3e+00
        &  \fat{1.3e+00}
        &  \fat{1.3e+00}
    \\ 
    & P-05
        & \fat{-3.7e+00}
        & \fat{-3.7e+00}
        & \fat{-3.7e+00}
        & -2.6e+00
        & -3.0e+00
        & -1.8e+00
        & -1.2e+00
        & \fat{-3.7e+00}
    \\ 
    & P-06
        &  \fat{1.2e-01}
        &  \fat{1.2e-01}
        &  \fat{1.2e-01}
        &  \fat{1.2e-01}
        &  \fat{1.2e-01}
        &  1.3e-01
        &  4.2e-01
        &  \fat{1.2e-01}
    \\ 
    & P-07
        &  \fat{6.2e+06}
        &  \fat{6.2e+06}
        &  \fat{6.2e+06}
        &  6.3e+06
        &  1.7e+07
        &  2.2e+10
        &  3.1e+08
        &  \fat{6.2e+06}
    \\ 
    & P-08
        &  \fat{6.0e+01}
        &  \fat{6.0e+01}
        &  \fat{6.0e+01}
        &  \fat{6.0e+01}
        &  \fat{6.0e+01}
        &  1.2e+02
        &  1.0e+02
        &  \fat{6.0e+01}
    \\ 
    & P-09
        &  \fat{2.7e+00}
        &  \fat{2.7e+00}
        &  \fat{2.7e+00}
        &  3.0e+00
        &  \fat{2.7e+00}
        &  2.9e+00
        &  3.4e+00
        &  \fat{2.7e+00}
    \\ 
    & P-10
        & \fat{-8.7e+02}
        & \fat{-8.7e+02}
        & \fat{-8.7e+02}
        & -6.1e+02
        & -6.9e+02
        &  7.0e+02
        &  2.6e+03
        & -8.5e+02
    \\ 

\hline
\multirow{4}{*}{QUBO}
    & P-11
        & \fat{-3.6e+02}
        & -3.5e+02
        & -3.4e+02
        & -3.2e+02
        & -3.4e+02
        & -3.2e+02
        & -3.3e+02
        & \fat{-3.6e+02}
    \\ 
    & P-12
        & \fat{-5.9e+03}
        & \fat{-5.9e+03}
        & \fat{-5.9e+03}
        & -5.6e+03
        & \fat{-5.9e+03}
        & -5.3e+03
        & \fat{-5.9e+03}
        & \fat{-5.9e+03}
    \\ 
    & P-13
        & \fat{-3.1e+00}
        & -3.0e+00
        & -2.8e+00
        &  0.0e+00
        &  1.5e+01
        &  2.8e+02
        & -2.9e+00
        & -3.0e+00
    \\ 
    & P-14
        & \fat{-3.1e+03}
        & -2.8e+03
        & -3.0e+03
        & -2.6e+03
        & -3.0e+03
        & -2.7e+03
        & -3.0e+03
        & -3.0e+03
    \\ 

\hline
\multirow{3}{*}{Control}
    & P-15
        &  \fat{6.7e-03}
        &  7.4e-03
        &  2.3e-02
        &  8.4e-03
        &  8.9e-03
        &  3.1e-02
        &  8.7e-02
        &  7.3e-03
    \\ 
    & P-16
        &  \fat{1.4e-02}
        &  2.6e-02
        &  3.5e-02
        &  1.7e-02
        &  1.7e-02
        &  5.3e-02
        &  5.2e-02
        &  \fat{1.4e-02}
    \\ 
    & P-17
        &  \fat{3.0e-02}
        &  5.7e-01
        &  1.5e-01
        &  4.8e-02
        &  3.6e-02
        &  7.7e-02
        &  5.3e-02
        &  3.7e-02
    \\ 

\hline
\multirow{3}{*}{\parbox{1.2cm}{Control +constr.}}
    & P-18
        &  1.4e-02
        &  \fat{1.1e-02}
        &  1.4e-02
        &  3.4e-02
        &  6.2e-02
        &  2.8e-01
        &  6.4e-02
        &  2.1e-02
    \\ 
    & P-19
        &  \fat{6.4e-02}
        &  5.7e-01
        &  6.7e-02
        & Fail
        & Fail
        & Fail
        & Fail
        & Fail
    \\ 
    & P-20
        &  \fat{1.5e-01}
        & Fail
        &  2.0e-01
        & Fail
        & Fail
        & Fail
        & Fail
        & Fail
    \\ \hline    
\end{tabular}
\end{sc}
\end{tiny}
\end{center}
\end{table}

\subsection{Multivariable analytic functions}
\label{sec:exp-analytic}

First, we consider the optimization task for various tensors
arising from the discretization of multivariable analytic functions.
We select $10$ popular benchmarks: Ackley (P-01), Alpine (P-02), Exponential (P-03), Griewank (P-04), Michalewicz (P-05), Piston\footnote{
    This function corresponds to the problem of modeling the time that takes a piston to complete one cycle within a cylinder; the description of its parameters can be found in \cite{zankin2018gradient, chertkov2022black}.

} (P-06), Qing (P-07), Rastrigin (P-08), Schaffer (P-09) and Schwefel (P-10).
These functions have a complex landscape and are often used in problems of evaluating the effectiveness of optimization algorithms~\cite{dieterich2012empirical,jamil2013literature}, including tensor-based optimizers~\cite{chertkov2022black, strossner2022approximation}.
We consider the 7-dimensional case (since this is the dimension of the Piston function) and discretization on a uniform grid with $16$ nodes.

As follows from Table~\ref{tbl:results} (benchmarks P-1, P-2, \ldots, P-10), our method, like the other two tensor approaches (BS-1 and BS-2), gave the most accurate solution for all model problems.
The most sophisticated approach from the nevergrad package (BS-7) turned out to be the next in accuracy (the method did not converge only in two cases out of ten).

\subsection{Quadratic unconstrained binary optimization}
\label{sec:exp-qubo}

QUBO is a widely known NP-hard problem~\cite{glover2022quantum} which unifies a rich variety of combinatorial optimization problems from finance and economics applications to machine learning and quantum computing.
QUBO formulation in a very natural manner utilizes penalty functions, yielding exact model representations in contrast to the approximate representations produced by customary uses of penalty functions.
The standard QUBO problem can be formulated as follows
\begin{equation}
    \ff(\vx) =
        \vx^T \matr{Q} \vx
        \rightarrow
        \min\limits_{\vx},
    \quad
    \textit{ s.t. }
    \vx \in \{ 0, 1\}^d,
\end{equation}
where $\vx$ is a vector of binary decision variables of the length $d$ and $\matr{Q} \in \set{R}^{d \times d}$ is a square matrix of constants.
In all our experiments we fixed the number of dimensions as $d = 50$.

We consider the following QUBO problems from the qubogen package:\footnote{
    See \url{https://github.com/tamuhey/qubogen}.
} Max-Cut Problem (P-11; which refers to finding a partition of an undirected graph into two sets such that the number of edges between the two sets is as large as possible), Minimum Vertex Cover Problem (P-12; which refers to finding a cover with a minimum number of vertices in the subset of the graph vertices such that each edge in the graph is incident) and Quadratic Knapsack Problem (P-13; which refers to finding a subset of maximum profit that satisfies the budget limitations from a set of potential projects with specified interactions between pairs of projects).

We also consider one more benchmark (P-14) from the work~\cite{dong2021phase} (problem $k_3$; $d = 50$), where angle-modulated bat algorithm (AMBA) was proposed for high-dimensional QUBO problems with engineering application to antenna topology optimization.
This is the ordinary binary knapsack problem with fixed weights $w_i \in [5, 20]$, profits $p_i \in [50, 100]$ ($i = 1, 2, \ldots, d$) and the maximum capacity $C = 1000$.
In experiments, we used the same values of the weights and profits as in~\cite{dong2021phase}.

For all four considered problems ( P-11, P-12, P-13, P-14) the proposed method PROTES gives the best result, as can be seen from Table~\ref{tbl:results}, and the baseline BS-7 again turned out to be the next in accuracy.
We also note that several optimization methods were compared in~\cite{dong2021phase} for the P-14 problem: BPSO (with the result $-2854$), BBA (with the result $-2976$), AMBA (with the result $-2956$), A-AMBA (with the result $-2961$), P-AMBA (with the result $-2989$), and the solution obtained using the PROTES method (the result $-3079$) turns out to be more accurate.

\subsection{Optimal control}
\label{sec:exp-control}



\def\iseq#1{#1(0), \,#1(1), \,\ldots, \,#1(T)}

Suppose we have a state variable $z \in \set{R}$ controlled by a binary variable $x$ called control (i.e., it's just a switch with modes ``off'' = 0 and ``on'' = 1) over some discrete interval of time $[0, T]$. 
The state~$z(t+1)$ at time $t+1$ 
depends on the control~$x(t)$ at time~$t$ and
obtained from the solution of the following differential equation $\dot{z}(\tau) = \ocf(z(\tau), x(t))$, $t \leq \tau < t+1$, where the function~$\ocf$ is called an equation function. 
The optimal control problem is to find such a sequence of controls $\vx^* = [\iseq{x^*}]$ (optimal solution) over the given time interval $[0, T]$ that minimizes the given objective function~$\oct$. 

Formulating the problem mathematically, we need to find such a solution
\begin{equation}
    \oct(\vz, \vx) \rightarrow
    \min\limits_{\vz, \vx},
    \quad
    \text{s.t.} \;\;
    \left\{ \begin{array}{l}
    z(0) = z_0, \\
    \dot z(\tau) = \ocf(z(\tau), \vx(t)),
        \;\; t\leq\tau< t+1,\\
    \vx(t) \in \{0, 1\},\;\; t=0,\,1,\,\ldots,\,T,
    \end{array} \right.
    \label{optimal-control}
\end{equation}
where $\vz = [z(0), z(1), \ldots, z(T)]$ is a state variable path.
In numerical experiments we consider the nonlinear equation function $\ocf(z, x) = z^3 - x$, and since it is nonlinear, finding an optimal solution raises a lot of difficulties.
The objective function $\oct$ we take in the form
$ 
    \oct(\vz, \vx) =
        \frac{1}{2}\sum\limits^{T}_{t=0}
            (z(t) - z_{\text{ref}})^2.
$ 
The initial and the reference state are fixed at values $z_0 = 0.8, \ z_{\text{ref}} = 0.7$.
For a fixed initial value~$z_0$ and fixed equation function~$\ocf$, the objective function~$\oct$ can be represented as a binary multidimensional tensor,
whose elements are calculated using the following function: $\ff(\vx) = \oct(\vz(\vx), \vx)$, hence we can apply discrete optimization methods to find $\vxmin$, which approximates the optimal solution $\vx^*$.

We considered several values for variable~$T$, 
such as $25, 50$ and $100$ (benchmarks P-15, P-16 and P-17 respectively).
As follows from the results presented in Table~\ref{tbl:results}, PROTES gives the most accurate solution in all three cases.
Note that a result comparable in accuracy for baselines is obtained only in one case.
Note that a result comparable in accuracy with our method is obtained only in one case when using baselines (i.e., P-16, BS-7).

\subsection{Optimal control with constraints}
\label{sec:exp-control_constr}

In practical applications, some conditions or constraints may be imposed on the solution of the optimal control problem.
We consider the following control constraint $P$ in this work: \textit{the control variable $\vx \in \{0, 1\}^N$ can take value ``1'' no less than $3$ times during the whole time interval}.
Formally, this can be written as follows:
\begin{equation*}
    P =
    \left\{
    \vx \ \bigg| \ 
    \begin{aligned}
        & x[t] \geq x[t-1] - x[t-2] \\
        & x[t] \geq x[t-1] - x[t-3]
    \end{aligned},
    \;\forall t:1\leq t\leq N+2;\text{\small we let $x[t]\coloneqq0$ for $t<1$ and $t>N$}
    \right\}.
\end{equation*}

To account for this condition in the PROTES algorithm,
we constructively build the initial distribution in the form of an indicator tensor as was described in Section~\ref{sec:scheme} in the constrained optimization subsection.
The details of this construction are presented in the Appendix.
The numerical results for $T=25$ (P-18), $T=50$ (P-19) and $T=100$ (P-20) are reported in Table~\ref{tbl:results}.
In two cases out of three (P-19, P-20), our method showed the best result, and in one case (P-18) slightly yielding to the TTOpt method (BS-1), which, however, in two other cases gave a significantly worse result.

\begin{figure}[t!]
    \centering
    \includegraphics[width=0.99\linewidth]
        {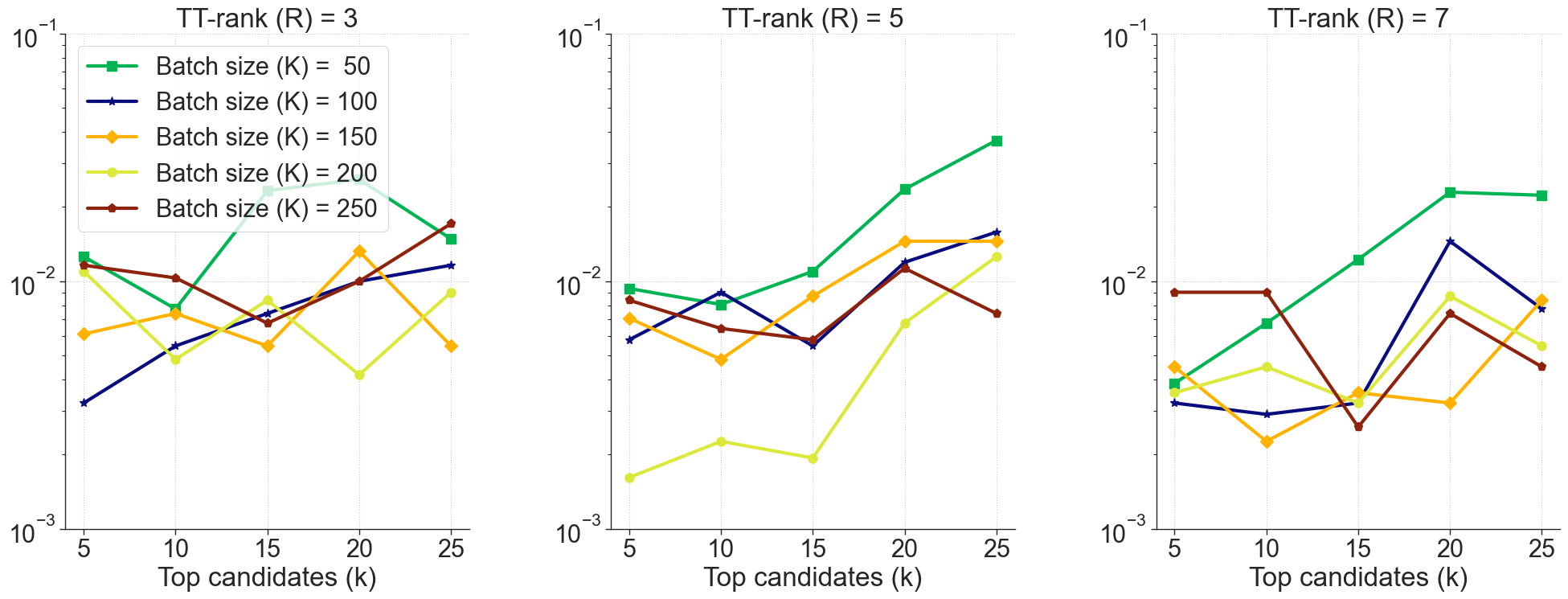}
    \caption{
        Relative error of the optimization result with PROTES method for the P-14 model problem with a known exact minimum for different values of the hyperparameters $K$, $k$ and $R$.
    }
    \label{fig:check}
\end{figure}

\subsection{Robustness and performance of the PROTES}
\label{sec:exp-performance}

The results in Table~\ref{tbl:results} relate to the ``intuitive'' selection of the hyperparameters for the PROTES method (as was mentioned above, we have used the values: $K=100$, $k=10$, $k_{gd}=1$, $\lambda=0.05$ and $R=5$).
In Figure~\ref{fig:check}, we present an analysis of the dependence of the optimization result for the benchmark P-14 on the choice of hyperparameters $K$, $k$ and $R$, with fixed $k_{gd}=1$ and $\lambda=0.05$.
We report the relative error of the result for all combinations $K = 50, 100, 150, 200, 250$; $k = 5, 10, 15, 20, 25$; $R = 3, 5, 7$.
As we can see from the plots, the hyperparameters used in the main calculations are not optimal for this particular problem, that is, additional fine-tuning of the method for specific problems or classes of problems is possible.
At the same time, according to the results in Figure~\ref{fig:check}, the method remains stable over a wide range of hyperparameter values.
We also note that all computations were carried out on a regular laptop, while the operating time of the considered optimizers was commensurate, for example, for the benchmark P-17, the measured operating time (in seconds) turned out to be as follows: PROTES (641), BS-1 (607), BS-2 (4245), BS-7 (780).
A more detailed analysis of the PROTES performance and the dependence of optimization results on the values of hyperparameters are considered in the Appendix.
\section{Conclusions}
    \label{sec:conclusions}
    In this work, we presented an optimization algorithm PROTES based on sampling from the probability density defined in the tensor train format.
For all considered numerical experiments, we used the same set of hyperparameters, so our algorithm is rather universal.
To take into account the constraints, as in the problem of optimal control with constraints, we only considered them in the form of a specially selected initial approximation (a special form of an indicator tensor in the tensor train format); further on, the algorithm did not consider the constraints explicitly.
This approach allows us to extend the capabilities of the algorithm by using the properties of the tensor train representation.
Numerical experiments show that we outperform many popular optimization methods.

The main direction in our future work is scaling of the method to large dimensions.
For $d \ge 1000$ we have encountered numerous technical difficulties, which can be alleviated by other tensor formats (such as hierachical Tucker, which can be parallelized over $d$) and more efficient implementations of the optimization method (now we used standard automatic differentiation without special tensor optimization methods such as Riemannian optimization).

\begin{ack}
    The work was supported by the Analytical center under the RF Government (subsidy agreement 000000D730321P5Q0002, Grant No. 70-2021-00145 02.11.2021).
\end{ack}
\printbibliography

@article{alarie2021two,
  title={Two decades of blackbox optimization applications},
  journal={EURO Journal on Computational Optimization},
  volume={9},
  pages={100011},
  year={2021},
  author={Alarie, Stephane and Audet, Charles and Gheribi, Aïmen and Kokkolaras, Michael and {Le Digabel}, Sebastien},
}

@article{alexandrov2023challenging,
  title={Challenging the curse of dimensionality in multidimensional numerical integration by using a low-rank tensor-train format},
  author={Alexandrov, Boian and Manzini, Gianmarco and Skau, Erik W and Truong, Phan Minh Duc and Vuchov, Radoslav G},
  journal={Mathematics},
  volume={11},
  number={3},
  pages={534},
  year={2023},
  publisher={Multidisciplinary Digital Publishing Institute}
}

@article{bennet2021nevergrad,
  author={Bennet, Pauline and Doerr, Carola and Moreau, Antoine and Rapin, Jeremy and Teytaud, Fabien and Teytaud, Olivier},
  title={Nevergrad: black-box optimization platform},
  year={2021},
  publisher={Association for Computing Machinery},
  volume={14},
  number={1},
  journal={SIGEVOlution},
  pages={8–15}
}

@article{chertkov2022black,
  title={Black box approximation in the tensor train format initialized by {ANOVA} decomposition},
  author={Chertkov, Andrei and Ryzhakov, Gleb and Oseledets, Ivan},
  journal={arXiv preprint arXiv:2208.03380},
  year={2022}
}

@article{chertkov2022optimization,
  title={Optimization of functions given in the tensor train format},
  author={Chertkov, Andrei and Ryzhakov, Gleb and Novikov, Georgii and Oseledets, Ivan},
  journal={arXiv preprint arXiv:2209.14808},
  year={2022}
}

@article{cichocki2016tensor,
  title={Tensor networks for dimensionality reduction and large-scale optimization: {Part} 1 low-rank tensor decompositions},
  author={Cichocki, Andrzej and Lee, Namgil and Oseledets, Ivan and Phan, Anh-Huy and Zhao, Qibin and Mandic, Danilo},
  journal={Foundations and Trends in Machine Learning},
  volume={9},
  number={4-5},
  pages={249--429},
  year={2016},
  publisher={Now Publishers Inc. Hanover, MA, USA}
}

@article{cichocki2017tensor,
  title={Tensor networks for dimensionality reduction and large-scale optimization: {Part} 2 applications and future perspectives},
  year={2017},
  volume={9},
  journal={Foundations and Trends in Machine Learning},
  number={6},
  pages={431-673},
  author={Andrzej Cichocki and Anh Phan and Qibin Zhao and Namgil Lee and Ivan Oseledets and Masashi Sugiyama and Danilo Mandic}
}

@article{dieterich2012empirical,
  title={Empirical review of standard benchmark functions using evolutionary global optimization},
  year={2012},
  journal={Applied Mathematics},
  author={Dieterich, Johannes and Hartke, Bernd},
  number={10},
  pages={1552--1564},
  volume={3}
}

@article{doerr2021survey,
  title={A survey on recent progress in the theory of evolutionary algorithms for discrete optimization},
  author={Doerr, Benjamin and Neumann, Frank},
  journal={ACM Transactions on Evolutionary Learning and Optimization},
  volume={1},
  number={4},
  pages={1--43},
  year={2021},
  publisher={ACM New York, NY}
}

@article{dolgov2020approximation,
  title={Approximation and sampling of multivariate probability distributions in the tensor train decomposition},
  author={Dolgov, Sergey and Anaya-Izquierdo, Karim and Fox, Colin and Scheichl, Robert},
  journal={Statistics and Computing},
  volume={30},
  pages={603--625},
  year={2020},
  publisher={Springer}
}

@article{dong2021phase,
  title={A phase angle-modulated bat algorithm with application to antenna topology optimization},
  author={Dong, Jian and Wang, Zhiyu and Mo, Jinjun},
  journal={Applied Sciences},
  volume={11},
  number={5},
  pages={2243},
  year={2021},
  publisher={MDPI}
}

@article{glover2022quantum,
  title={Quantum bridge analytics I: a tutorial on formulating and using {QUBO} models},
  author={Glover, Fred and Kochenberger, Gary and Hennig, Rick and Du, Yu},
  journal={Annals of Operations Research},
  pages={1--43},
  year={2022},
  publisher={Springer}
}

@article{grasedyck2013literature,
  author={Grasedyck, Lars and Kressner, Daniel and Tobler, Christine},
  title={A literature survey of low-rank tensor approximation techniques},
  journal={GAMM-Mitteilungen},
  volume={36},
  number={1},
  pages={53--78},
  year={2013}
}

@article{hammersley1960monte,
  title={Monte {C}arlo methods for solving multivariable problems},
  author={Hammersley, John},
  journal={Annals of the New York Academy of Sciences},
  volume={86},
  number={3},
  pages={844--874},
  year={1960},
  publisher={Wiley Online Library}
}

@inproceedings{hansen2006cma,
  author={Hansen, Nikolaus},
  title={The CMA Evolution Strategy: A Comparing Review},
  booktitle={Towards a new evolutionary computation: advances in the estimation of distribution algorithms},
  year={2006},
  publisher={Springer Berlin Heidelberg},
  pages={75--102}
}

@article{jamil2013literature,
  title={A literature survey of benchmark functions for global optimization problems},
  year={2013},
  journal={Journal of Mathematical Modelling and Numerical Optimisation},
  author={Jamil, Momin and Yang, Xin-She},
  number={2},
  pages={150--194},
  volume={4}
}

@inproceedings{kennedy1995particle,
  title={Particle swarm optimization},
  author={Kennedy, James and Eberhart, Russell},
  booktitle={Proceedings of ICNN'95-international conference on neural networks},
  volume={4},
  pages={1942--1948},
  year={1995},
  organization={IEEE}
}

@article{kingma2014adam,
  title={Adam: a method for stochastic optimization},
  author={Kingma, Diederik P and Ba, Jimmy},
  journal={arXiv preprint arXiv:1412.6980},
  year={2014}
}

@article{kour2023efficient,
  title={Efficient structure-preserving support tensor train machine},
  author={Kour, Kirandeep and Dolgov, Sergey and Stoll, Martin and Benner, Peter},
  journal={Journal of Machine Learning Research},
  volume={24},
  number={4},
  pages={1--22},
  year={2023}
}

@inproceedings{lauand2022approaching,
  title={Approaching quartic convergence rates for quasi-stochastic approximation with application to gradient-free optimization},
  author={Lauand, Caio Kalil and Meyn, Sean P},
  booktitle={Advances in Neural Information Processing Systems},
  year={2022}
}

@article{lin2022gradientfree,
  title={Gradient-free methods for deterministic and stochastic nonsmooth nonconvex optimization},
  author={Lin, Tianyi and Zheng, Zeyu and Jordan, Michael},
  journal={Advances in Neural Information Processing Systems},
  volume={35},
  pages={26160--26175},
  year={2022}
}

@inproceedings{maryak2001global,
  title={Global random optimization by simultaneous perturbation stochastic approximation},
  author={Maryak, John and Chin, Daniel},
  booktitle={Proceedings of the 2001 American Control Conference.(Cat. No. 01CH37148)},
  volume={2},
  pages={756--762},
  year={2001},
  organization={IEEE}
}

@article{nesterov2017random,
  author={Nesterov, Yurii and Spokoiny, Vladimir},
  title={Random gradient-free minimization of convex functions},
  journal={Foundations of Computational Mathematics},
  year={2017},
  volume={17},
  number={2},
  pages={527-566}
}

@article{nikitin2022quantum,
  title={Are quantum computers practical yet? {A} case for feature selection in recommender systems using tensor networks},
  author={Nikitin, Artyom and Chertkov, Andrei and Ballester-Ripoll, Rafael and Oseledets, Ivan and Frolov, Evgeny},
  journal={arXiv preprint arXiv:2205.04490},
  year={2022}
}

@inproceedings{novikov2021tensor,
  title={Tensor-train density estimation},
  author={Novikov, Georgii and Panov, Maxim and Oseledets, Ivan},
  booktitle={Uncertainty in artificial intelligence},
  pages={1321--1331},
  year={2021},
  organization={PMLR}
}

@article{oseledets2009breaking,
  title={Breaking the curse of dimensionality, or how to use {SVD} in many dimensions},
  author={Oseledets, Ivan and Tyrtyshnikov, Eugene},
  journal={SIAM Journal on Scientific Computing},
  volume={31},
  number={5},
  pages={3744--3759},
  year={2009},
  publisher={SIAM}
}

@article{oseledets2010ttcross,
  title={{TT}-cross approximation for multidimensional arrays},
  author={Oseledets, Ivan and Tyrtyshnikov, Eugene},
  journal={Linear Algebra and its Applications},
  volume={432},
  number={1},
  pages={70--88},
  year={2010},
  publisher={Elsevier}
}

@article{oseledets2011tensor,
  title={Tensor-train decomposition},
  author={Oseledets, Ivan},
  journal={SIAM Journal on Scientific Computing},
  volume={33},
  number={5},
  pages={2295--2317},
  year={2011},
  publisher={SIAM}
}

@article{parcollet2023learning,
  title={Learning {F}eynman diagrams with tensor trains},
  author={Parcollet, Olivier and Nunez-Fernandez, Yuriel and Jeannin, Matthieu and Dumitrescu, Philipp and Kloss, Thomas and Kaye, Jason and Waintal, Xavier},
  journal={Bulletin of the American Physical Society},
  year={2023},
  publisher={APS}
}

@book{parker2014discrete,
  title={Discrete optimization},
  author={Parker, Gary and Rardin, Ronald},
  year={2014},
  publisher={Elsevier}
}

@article{qi2023exploiting,
  title={Exploiting low-rank tensor-train deep neural networks based on {R}iemannian gradient descent with illustrations of speech processing},
  author={Qi, Jun and Yang, Chao-Han Huck and Chen, Pin-Yu and Tejedor, Javier},
  journal={IEEE/ACM Transactions on Audio, Speech, and Language Processing},
  volume={31},
  pages={633--642},
  year={2023},
  publisher={IEEE}
}

@inproceedings{ryzhakov2023constructive,
  title={Constructive {TT}-representation of the tensors given as index interaction functions with applications},
  author={Ryzhakov, Gleb and Oseledets, Ivan},
  booktitle={11th International Conference on Learning Representations, {ICLR}},
  year={2023}
}

@inproceedings{scarlett2017lower,
  title={Lower bounds on regret for noisy gaussian process bandit optimization},
  author={Scarlett, Jonathan and Bogunovic, Ilija and Cevher, Volkan},
  booktitle={Conference on Learning Theory},
  pages={1723--1742},
  year={2017},
  organization={PMLR}
}

@article{selvanayagam2022global,
  title={Global optimization of surface warpage for inverse design of ultra-thin electronic packages using tensor train decomposition},
  author={Selvanayagam, Cheryl and Duong, Pham Luu Trung and Wilkerson, Brett and Raghavan, Nagarajan},
  journal={IEEE Access},
  volume={10},
  pages={48589--48602},
  year={2022},
  publisher={IEEE}
}

@article{shetty2022tensor,
  title={Tensor train for global optimization problems in robotics},
  author={Shetty, Suhan and Lembono, Teguh and Loew, Tobias and Calinon, Sylvain},
  journal={arXiv preprint arXiv:2206.05077},
  year={2022}
}

@article{soley2023global,
  title={Global optimization with the iterative power algorithm via quantum computing and quantics tensor trains},
  author={Soley, Micheline and Kyaw, Thi Ha and Bergold, Paul and Allen, Brandon and Sun, Chong and Aspuru-Guzik, Al{\'a}n and Batista, Victor},
  journal={Bulletin of the American Physical Society},
  year={2023},
  publisher={APS}
}

@inproceedings{sozykin2022ttopt,
  title={{TTOpt}: a maximum volume quantized tensor train-based optimization and its application to reinforcement learning},
  author={Sozykin, Konstantin and Chertkov, Andrei and Schutski, Roman and Phan, Anh-Huy and Cichocki, Andrzej and Oseledets, Ivan},
  booktitle={Advances in Neural Information Processing Systems},
  year={2022}
}

@article{storn1997differential,
  author={Storn, Rainer and Price, Kenneth},
  title={Differential evolution -- a simple and efficient heuristic for global optimization over continuous spaces},
  journal={Journal of Global Optimization},
  year={1997},
  volume={11},
  number={4},
  pages={341-359}
}

@article{strossner2022approximation,
  title={Approximation in the extended functional tensor train format},
  author={Str{\"o}ssner, Christoph and Sun, Bonan and Kressner, Daniel},
  journal={arXiv preprint arXiv:2211.11338},
  year={2022}
}

@article{wang2022zarts,
  title={{ZARTS}: on zero-order optimization for neural architecture search},
  author={Wang, Xiaoxing and Guo, Wenxuan and Su, Jianlin and Yang, Xiaokang and Yan, Junchi},
  journal={Advances in Neural Information Processing Systems},
  volume={35},
  pages={12868--12880},
  year={2022}
}

@article{wang2023tensor,
  title={Tensor decompositions for hyperspectral data processing in remote sensing: a comprehensive review},
  author={Wang, Minghua and Hong, Danfeng and Han, Zhu and Li, Jiaxin and Yao, Jing and Gao, Lianru and Zhang, Bing and Chanussot, Jocelyn},
  journal={IEEE Geoscience and Remote Sensing Magazine},
  year={2023},
  publisher={IEEE}
}

@article{wierstra2014natural,
  author={Wierstra, Daan and Schaul, Tom and Glasmachers, Tobias and Sun, Yi and Peters, Jan and Schmidhuber, J\"{u}rgen},
  title={Natural evolution strategies},
  journal={Journal of Machine Learning Research},
  year={2014},
  volume={15},
  number={27},
  pages={949-980}
}

@article{williams1992simple,
  title={Simple statistical gradient-following algorithms for connectionist reinforcement learning},
  author={Williams, Ronald},
  journal={Machine learning},
  volume={8},
  number={3-4},
  pages={229--256},
  year={1992},
  publisher={Springer}
}

@article{zankin2018gradient,
  title={Gradient descent-based {D}-optimal design for the least-squares polynomial approximation},
  author={Zankin, Vitaly and Ryzhakov, Gleb and Oseledets, Ivan},
  journal={arXiv preprint arXiv:1806.06631},
  year={2018}
}
\vspace{1.5cm}
\centerline{\textbf{\large{Supplementary Material}}}
\setcounter{section}{0}

\section{Stability of the PROTES}
\label{s:appendix_stability}

To check the stability of the optimization result, a series of $10$ calculations were performed for each method (PROTES, BS1 -- BS7) with random initializations.
We consider the binary knapsack problem from the work~\cite{dong2021phase} (benchmark P-14, described in detail in the main text of the work) with the known exact minimum $-3103$.
For all methods, the limit on the number of requests was fixed at the value $M = 10^5$.
All other parameters were the same as in the computations from the main text, i.e., the PROTES parameters are $K = 100$, $k = 10$, $k_{gd} = 1$, $\lambda = 0.05$, $R = 5$.
In Table~\ref{tbl:results_seed} we present the average and best results over $10$ runs for each optimization method.
As follows from the reported results, only PROTES and Portfolio method (BS-7) managed to successfully find the exact optimum, while the average result for PROTES is significantly better than that of all the baselines.

\begin{table}[h!]
\caption{
    Average and best result for $10$ independent runs for the P-14 benchmark.
}
\label{tbl:results_seed}
\renewcommand{\arraystretch}{1.5}
\begin{center}
\begin{small}
\begin{sc}
\begin{tabular}
{p{1.0cm}rrrrrrrr}\hline

        &
PROTES  & 
BS-1    & 
BS-2    &
BS-3    &
BS-4    &
BS-5    &
BS-6    &
BS-7 
\\ \hline
Mean
        & \fat{-3095}
        & -2992
        & -3048
        & -2650
        & -2937
        & -2701
        & -3064
        & -3075
\\
Best
        & \fat{-3103}
        & -3074
        & -3084
        & -2825
        & -2996
        & -2752
        & -3094
        & \fat{-3103}
\\  \hline
\end{tabular}
\end{sc}
\end{small}
\end{center}
\end{table}

\section{Choice of hyperparameters for the PROTES}
\label{s:appendix_hyperparameters}

\begin{figure}[t!]
    \centering
    \includegraphics
        [scale=0.27]
        {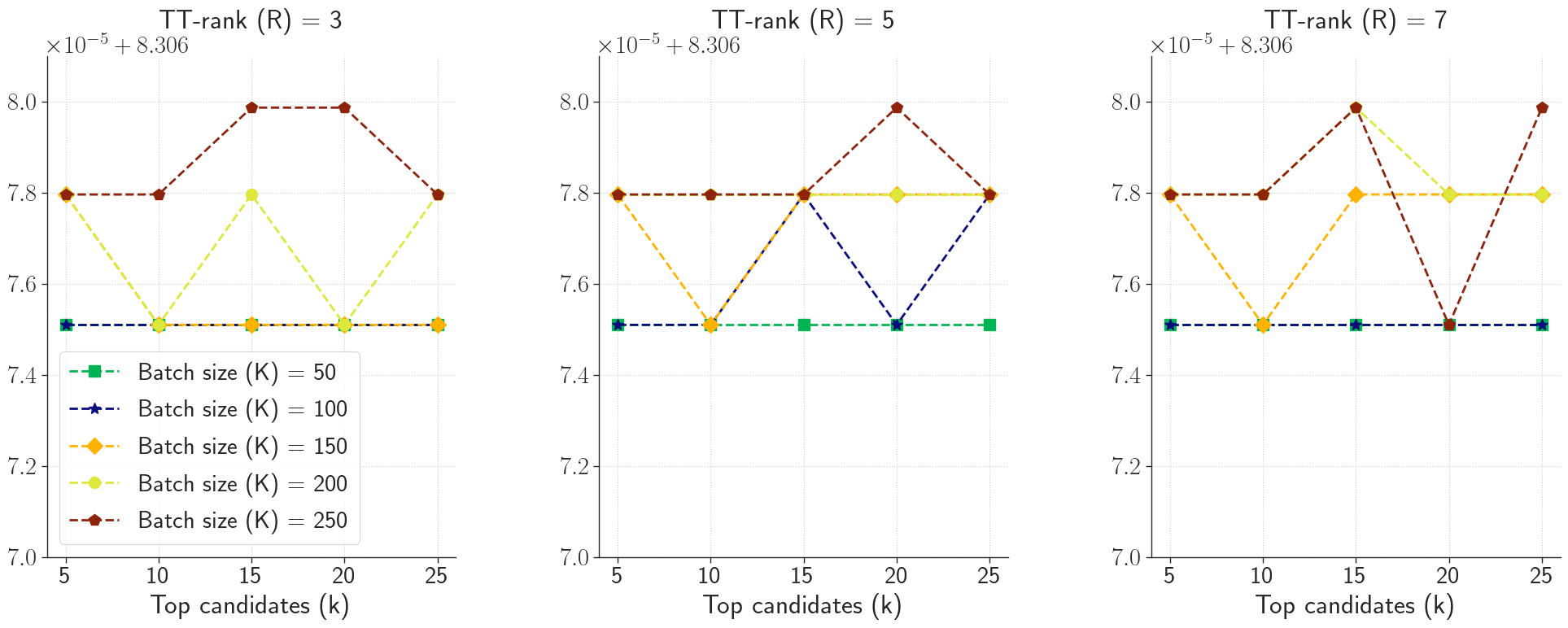}
    \caption{
        Optimization result (approximation of the minimum value) with PROTES method for the P-01 for different values of the hyperparameters: $K$, $k$ and $R$.
    }
    \label{fig:check_ackley}
\end{figure}

\begin{figure}[t!]
    \centering
    \includegraphics
        [scale=0.27]
        {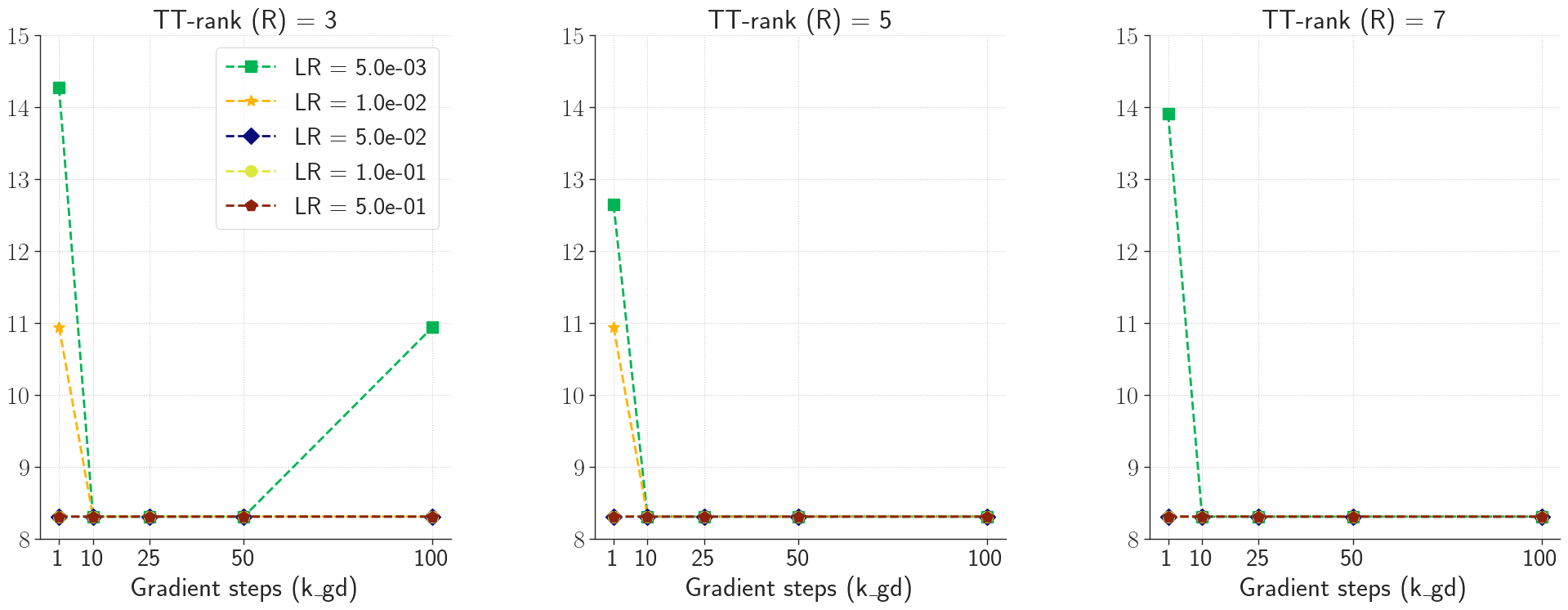}
    \caption{
        Optimization result (approximation of the minimum value) with PROTES method for the P-01 model problem for different values of the hyperparameters: learning rate (LR), $k_{gd}$, $R$.
    }
    \label{fig:check_ackley_ext}
\end{figure}

\begin{figure}[t!]
    \centering
    \includegraphics
        [scale=0.27]
        {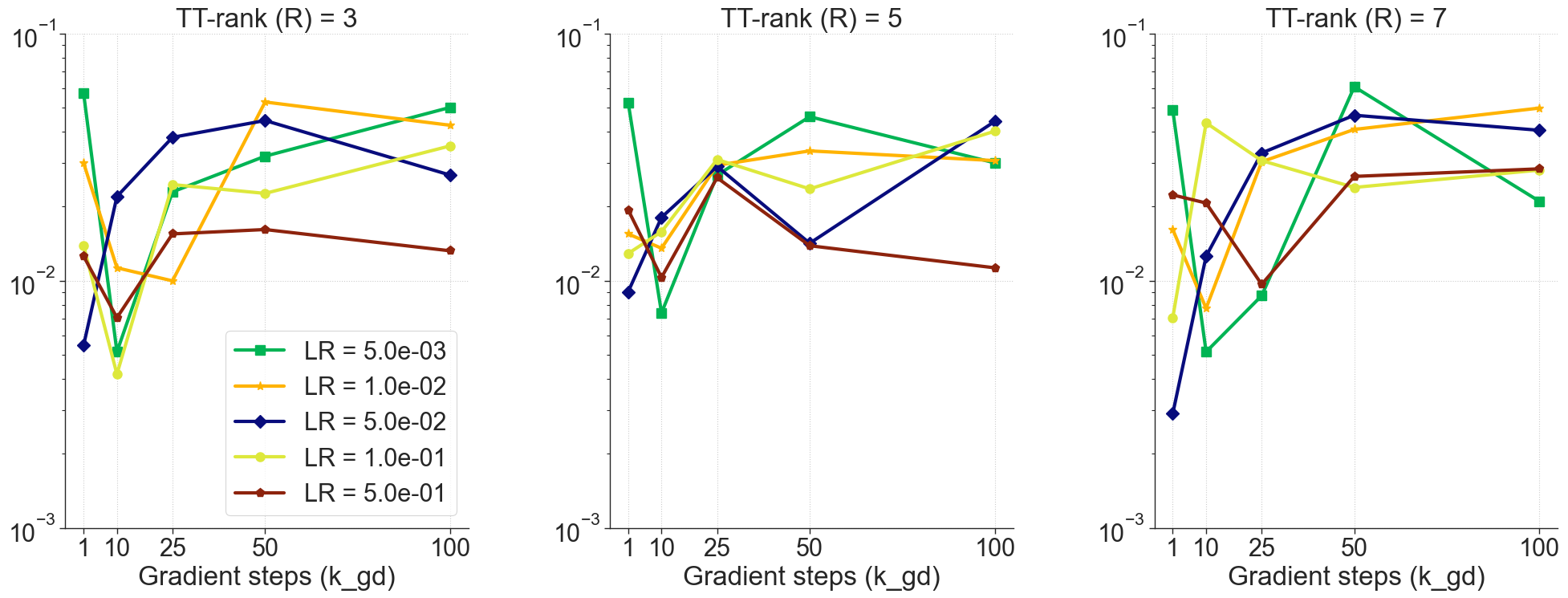}
    \caption{
        Relative error of the optimization result with PROTES method for the P-14 model problem with a known exact minimum for different values of the hyperparameters: learning rate (LR), $k_{gd}$, $R$.
    }
    \label{fig:check_knapsack_ext}
\end{figure}

We conduct additional experiments with varying the values of the hyperparameters of the PROTES method for benchmarks P-01 (Figure~\ref{fig:check_ackley} and Figure~\ref{fig:check_ackley_ext}; we report the optimization result) and P-14 (Figure~4 in the main text and Figure~\ref{fig:check_knapsack_ext} below; we report the relative error of the optimization result).
In the first series of experiments (Figure~\ref{fig:check_ackley} below and Figure~4 in the main text), we fixed the values $k_{gd} = 1$ and $\lambda = 0.05$ and tried all combinations of the remaining hyperparameters: $K = 50, 100, 150, 200, 250$; $k = 5, 10, 15, 20, 25$; $R = 3, 5, 7$.
In the second series of experiments (Figure~\ref{fig:check_ackley_ext} and~\ref{fig:check_knapsack_ext}), we fixed the values $K = 100$ and $k = 10$, and tried all combinations: $k_{gd} = 1, 10, 25, 50, 100$, $\lambda = 0.005, 0.01, 0.05, 0.1, 0.5$; $R = 3, 5, 7$.

As follows from the presented results, for problem P-01, the dependence on the choice of hyperparameters turns out to be extremely weak, except for outliers for the small values of the learning rate $\lambda = 0.005$ and $\lambda = 0.01$.
For more complex problem P-14, the dependence of the result on the choice of hyperparameters is more complex, but it can be seen that the method remains stable over a wide range of hyperparameter values.

\section{Performance comparison of optimization methods}
\label{s:appendix_performance}

All the results described in the main text and presented there in Table~1 were obtained on a regular laptop.
In Table~\ref{tbl:results_time} we report the related computation time for each method (PROTES, BS1 -- BS7) and each model problem (P01 -- P20).
As follows from the results, the PROTES works much faster than classical optimization methods (BS3 -- BS7), as well as faster than tensor-based methods (BS1 and BS2) for most optimal control problems (P-15 -- P-20).
However, the TTOpt (BS-1) and Optima-TT (BS-2) methods are faster for simpler analytic (P-01 -- P-10) and QUBO (P-11 -- P-14) problems.
This is due to the fact that within the framework of the dynamic TT-rank refinement procedure in these methods, the TT-rank, and hence the computation time, turn out to be significantly higher for the optimal control problems.
We note that the very short running time of the TTOpt method for P-20 is because most of the requests of the method did not satisfy the imposed constraints, and in this case, the differential equation was not solved.

\begin{table}[t!]
\caption{
    Computation time in seconds for all selected benchmarks (P-01 -- P-20) and for all used optimization methods (PROTES, BS1 -- BS7).
}
\label{tbl:results_time}
\renewcommand{\arraystretch}{1.3}
\begin{center}
\begin{tiny}
\begin{sc}
\begin{tabular}
{p{2.0cm}p{1.0cm}rrrrrrrr}\hline

        &
        &
PROTES  & 
BS-1    & 
BS-2    &
BS-3    &
BS-4    &
BS-5    &
BS-6    &
BS-7    \\ \hline

\multirow{10}{*}{\parbox{1.3cm}{Analytic Functions}}
    & P-01
        & 3.28
        & 0.06
        & 0.11
        & 23.45
        & 25.31
        & 22.78
        & 22.03
        & 67.23
    \\
    & P-02
        & 2.25
        & 0.05
        & 0.06
        & 22.42
        & 28.68
        & 20.6
        & 19.78
        & 77.01
    \\
    & P-03
        & 2.36
        & 0.05
        & 0.03
        & 26.11
        & 23.62
        & 20.55
        & 19.71
        & 65.62
    \\
    & P-04
        & 2.32
        & 0.05
        & 0.07
        & 22.04
        & 24.04
        & 20.67
        & 19.9
        & 65.51
    \\
    & P-05
        & 2.33
        & 0.05
        & 0.06
        & 26.47
        & 26.56
        & 20.78
        & 20.1
        & 78.42
    \\
    & P-06
        & 2.34
        & 0.05
        & 0.1
        & 23.51
        & 28.32
        & 22.42
        & 21.58
        & 84.1
    \\
    & P-07
        & 2.25
        & 0.06
        & 0.03
        & 26.72
        & 34.95
        & 21.66
        & 21.26
        & 83.54
    \\
    & P-08
        & 2.25
        & 0.05
        & 0.07
        & 22.7
        & 24.03
        & 21.29
        & 20.28
        & 65.47
    \\
    & P-09
        & 2.31
        & 0.06
        & 0.11
        & 22.98
        & 24.11
        & 21.14
        & 20.56
        & 65.55
    \\
    & P-10
        & 2.35
        & 0.05
        & 0.03
        & 23.56
        & 33.66
        & 21.25
        & 22.61
        & 81.94
    \\

\hline
\multirow{4}{*}{QUBO}
    & P-11
        & 2.7
        & 0.06
        & 0.44
        & 16.32
        & 22.54
        & 17.65
        & 17.98
        & 67.61
    \\
    & P-12
        & 2.16
        & 0.05
        & 0.36
        & 17.06
        & 21.95
        & 17.28
        & 16.77
        & 79.74
    \\
    & P-13
        & 2.29
        & 0.05
        & 0.34
        & 18.78
        & 21.09
        & 20.51
        & 17.97
        & 77.22
    \\
    & P-14
        & 2.33
        & 0.13
        & 0.4
        & 16.41
        & 27.42
        & 17.61
        & 17.59
        & 74.77
    \\

\hline
\multirow{3}{*}{Control}
    & P-15
        & 513.6
        & 1256.0
        & 1839.0
        & 550.0
        & 545.0
        & 620.3
        & 530.4
        & 707.5
    \\
    & P-16
        & 542.4
        & 969.4
        & 3007.0
        & 578.0
        & 595.0
        & 570.7
        & 573.7
        & 697.5
    \\
    & P-17
        & 640.7
        & 607.3
        & 4245.0
        & 673.6
        & 687.6
        & 661.0
        & 687.1
        & 779.5
    \\

\hline
\multirow{3}{*}{\parbox{1.2cm}{Control +constr.}}
    & P-18
        & 328.1
        & 92.9
        & 588.1
        & 319.1
        & 516.2
        & 202.9
        & 533.9
        & 474.5
    \\
    & P-19
        & 8.74
        & 69.49
        & 912.0
        & 17.13
        & 17.45
        & 16.16
        & 16.84
        & 47.43
    \\
    & P-20
        & 9.23
        & 0.53
        & 931.8
        & 21.16
        & 21.42
        & 20.28
        & 20.84
        & 62.26
    \\  \hline
\end{tabular}
\end{sc}
\end{tiny}
\end{center}
\end{table}

\section{Derivative functions for the indicator tensor in the optimal control problem}
\label{s:appendix_derivative}

We use the following derivative functions for the constructive building of the tensor described in~\cite{ryzhakov2023constructive}
\begin{equation}
\begin{split}
f^k_0(x) & =
    \left\{
        \begin{aligned}
            0&,   &x=0 \text{ or } x=l\\
            \texttt{None}&, &\text{otherwise},
        \end{aligned}
    \right.
\\
f^k_1(x) & = \min(l,\, x+1),
\end{split}
\end{equation}
for all TT-cores except the last one ($k = 1, 2, \ldots, d-1$), and for the last TT-core
\begin{equation}
\begin{split}
f^d_0(x) & =
    \left\{
        \begin{aligned}
            1&,   &x=0 \text{ or } x=l\\
            0&, &\text{otherwise},
        \end{aligned}
    \right.
\\
f^d_1(x) & =
    \left\{
        \begin{aligned}
            1&, &x \geq l-1\\
            0&, &\text{otherwise}.
        \end{aligned}
    \right.
\end{split}
\end{equation}
A tensor in the TT-format built on such derivative functions is equal to $0$ if there are less than~$l$ ones among its vector argument in a row, and is equal to $1$ in all other cases.

Let us briefly explain why such derivative functions give a tensor of the restriction condition (i.e., the constraint in the considered optimal control problem).
Recall that the upper index in the derivative function notation corresponds to the index number of the tensor argument, and the lower index corresponds to the value of this argument index.
In our scheme, the argument of the derivative functions has the meaning of the number of ones that already stand to the left of the given index.

First, let's focus on the functions for all indices except the rightmost one ($k = 1, 2, \ldots, d-1$).
If the current index is one, then we simply increase the value of the argument by~$1$ and pass it on, to the input of the next derivative function.
This is what the function~$f^k_1$ does (we take the maximum, so if the number of ones has already reached~$l$, we don't care if it is greater than or equal to~$l$; note that the maximum operation reduces the TT-rank but does not affect the result).
For the zero value of the current index, if the argument is zero (which means, that the previous index is also zero) the function~$f^k_0$ also returns zero as nothing have changed.
If the argument of this function is~$l$, it means that there are~$l$ or more ones in a row to the left of the considered index, so the given condition is not violated, and the function simply returns zero, which means that there are no ones.
If the argument is greater than zero but less than~$l$, it means that the condition is violated, because the previous sequence of ones of length less than $l$ is cut off at the current index. 
In this case, the function~$f^k_0$ returns \texttt{None}, which means that the value of the tensor will be~$0$ regardless of the subsequent indices.

The functions that correspond to the last ($k=d$) index behave similarly.
Namely, if the last index is zero, then the function~$f^k_0$ returns~$1$ if its argument is~$0$ or~$l$ since it does not violate the condition as just explained.
Otherwise, it returns~$0$.
If the last index is~$1$, then the function~$f^k_1$ returns~$1$ only if its argument is equal to $(l-1)$ or~$l$ since this means that the current one is an element of a sequence of ones of length at least~$l$.
Otherwise, the function~$f^k_1$ returns~$0$.
\end{document}